\documentstyle[12pt]{article} \vsize=28.7truecm \hsize=23truecm
 \def\vbar{\mathchoice{\vrule height2.3ptdepth-.3ptwidth.12pt\kern-
 .10pt}
    {\vrule height6.3ptdepth-.3ptwidth.11pt\kern-.11pt}
    {\vrule height5.1ptdepth-.30ptwidth.8pt\kern-.8pt}
    {\vrule height4.1ptdepth-.24ptwidth.6pt\kern-.7pt}}
\def\reel{\hbox{I\hskip-2pt R}}

\def\<{\langle}
\def\>{\rangle}

\def\n{{\boldmath n}}

\def\<{\langle}
\def\>{\rangle}
\def\reel{\hbox{I\hskip -2pt R}}

\def\n{{\noindent}}
\newtheorem{theorem}{Theorem}[section]

\newtheorem{prop}[theorem]{Proposition}

\newenvironment{Pff}{\hspace*{-\parindent}{\bf Proof}}{\hfill $\Box$
\vspace*{0.2cm}}

\begin{document}
 \vspace{2cm}

\title{{\bf{New characterizations of the Letac-Mora class of real cubic natural exponential
families}}}
\author{Marwa Hamza and Abdelhamid Hassairi \footnote{Corresponding author.
 \textit{E-mail address: Abdelhamid.Hassairi@fss.rnu.tn}}$\;$
\\{\footnotesize{\it
Laboratory of Probability and Statistics. Sfax Faculty of Sciences,
B.P. 1171, Tunisia.}}}
\date{}
\maketitle \n $\overline{\hspace{15cm}}$\vskip0.3cm \n {\small {\bf
Abstract}} {\small In this paper, we give three equivalent
properties of the class of multivariate simple cubic natural
exponential families (NEF's). The first property says that the
cumulant function of any basis of the family is a solution of some
Monge-Amp\`{e}re equation, the second is that the variance function
satisfies a differential equation, and the third is characterized by
the equality between two families of prior distributions related to
the NEF. These properties represent the extensions to this class of
the properties stated in $\cite{Casalis(1996)}$ and satisfied by the
Wishart and the simple quadratic NEF's. We also show that in the
real case, each of these properties provides a new characterization
of the Letac-Mora class of real cubic NEF's.}
\vskip0.2cm\n{\small{\it{Keywords: Natural exponential family,
variance function, cumulant function, Monge-Amp\`{e}re equation, prior distribution.}} \\
$\overline{\hspace{15cm}}$\vskip1cm

\section{Introduction and preliminaries}

For the convenance of the reader, we first introduce some notations
and recall some facts concerning the natural exponential families
and their variance functions, our notations are the ones used in
$\cite{G.}$. Let $E$ be a linear vector space with finite dimension
$n,$ denote by $E^{*}$ its dual, and let $E^{*}\times E \rightarrow
\reel:(\theta,x)\mapsto
\langle\theta,x\rangle$ be the duality bracket.  \\
If $\mu$ is a positive radon measure on $E$, then
\begin{equation}\label{MI1}
L_\mu(\theta)=\displaystyle\int_{E} \exp(\langle\theta,
x\rangle)\mu(dx)
\end{equation}
denotes its Laplace transform. We also denote by ${\mathcal{M}}(E)$
the set of measures $\mu$ such that the set
\begin{equation}
\Theta(\mu)=\textrm{interior}\{\theta \in E^{*}; \ L_\mu
(\theta)<+\infty\}
\end{equation}
is non empty and $\mu$ is not concentrated on an affine hyperplane
of $E$. The cumulant function of an element $\mu$ of
${\mathcal{M}}(E)$ is the function defined for $\theta \textrm{ in
}\Theta(\mu)$ by
$$ k_\mu (\theta)=\log L_\mu (\theta).$$ \vskip0.1cm\n To
each $\mu$ in ${\mathcal{M}}(E)$ and $\theta$ in $\Theta(\mu)$, we
associate the probability distribution on $E$ \vskip0.2cm$\hfill
P(\theta,\mu)(dx)=exp\left(\langle\theta,x\rangle -k_\mu
(\theta)\right)\mu(dx).\hfill$\newline
 The set\vskip0.1cm$\hfill
 F=F(\mu)=\{P(\theta,\mu); \ \theta\in\Theta(\mu)\}\hfill$
\vskip0.3cm\n is called the natural exponential family (NEF)
generated by $\mu$. We also say that $\mu$ is a basis of $F$.
\vskip0.2cm\n The function $k_\mu$ is strictly convex and real
analytic. Its first derivative $k'_\mu$ defines a diffeomorphism
between $\Theta(\mu)$ and its image $M_F$. Since $k'_\mu
(\theta)=\displaystyle\int_E x P(\theta,\mu) (dx)$, $M_F$ is called
the domain of the means of $F$. The inverse function of $k'_\mu$ is
denoted by $\psi_\mu$ and setting $$P(m,F)=P(\psi_\mu(m),\mu)$$ the
probability of $F$ with mean $m$, we have $$F=\left\{P(m,F); \ m\in
M_F\right\},$$ which is the parametrization of $F$ by the
mean.\vskip0.1cm\n Now the covariance operator of $P(m,F)$ is
denoted by $V_F (m)$ and the map $$M_F\longrightarrow L_s(E^{*},E);
\ m\longmapsto V_F (m)=k''_\mu(\psi_\mu(m))$$ is called the variance
function of the NEF $F$. It is easy proved that for all $m\in M_F$,
$$V_F(m)=(\psi'_\mu(m))^{-1},$$ and an important feature of $V_F$ is
that it characterizes $F$ in the following sense: If $F$ and $F'$
are two NEFs such that $V_F (m)$ and $V_{F'}(m)$ coincide on a
nonempty open set of $M_F \cap M_{F'},$ then $F=F'$. \\ Now, let us
examine the influence of an affine transformation and a power
convolution on a NEF $F=F(\mu)$. If $\varphi(x)=a(x)+b$, where $a\in
GL(E)$ and $b\in E$, is an affine transformation of $E$, then
$\varphi(F(\mu))=F(\varphi(\mu)),$ $M_{\varphi(F)}=\varphi(M_{F}),$
and $$V_{\varphi(F)}(m)=a\ V_{F}(\varphi^{-1}(m))\ a^{*},$$ where
$a^{*}$ is the transpose of $a$. On the other hand the set
$$
\Lambda(\mu)=\{\lambda>0; \ \exists \ \mu_{\lambda} \in
{\mathcal{M}}(E) \textrm{ such that }
L_{\mu_{_{\lambda}}}(\theta)=\left(L_{\mu}(\theta)\right)^{\lambda}
\textrm{ for all }\theta \in \Theta(\mu) \}$$ is called the
Jorgensen set of $\mu$ and the measure $\mu_{\lambda}$ is the
$\lambda-$power of convolution of $\mu$. For $\lambda$ in
$\Lambda(\mu)$, we have that
$$M_{F(\mu_{\lambda})}=\lambda M_{F}, \textrm{ and }
V_{F_{\lambda}}(m)=\lambda\ V_{F}(\displaystyle\frac{m}{\lambda}).$$

A very interesting fact is that the most common real and
multivariate probability distributions belong to the natural
exponential families such that the variance function is a polynomial
of degree less then or equal to three in the mean $m$. For instance,
up to affine transformations and power of convolution (up to the
type), the Gaussian, Poisson, gamma, binomial, negative binomial and
hyperbolic cosine distributions form the class of all real NEF's
whose variance function is a polynomial of degree less than or equal
to 2 characterized by Morris $\cite{Morris(1982)}$. Letac and Mora
$\cite{Letac and Mora (1990)}$ have added six others types of
distributions, namely, the inverse Gaussian, Ressel, Abel,
Tack$\grave{a}$s, strict arcsine and large arcsine, to get the class
of real cubic NEF's, that is the class of NEF's such that variance
function is a polynomial of degree less than or equal to three. The
classification of NEF's with polynomial variance function have been
extended to the multivariate NEF's. The multivariate version of the
Morris class, called the class of simple quadratic NEF's, has been
completely described by Casalis $\cite{Casalis(1996)}$, it contains
$2n+4$ types. Hassairi $\cite{Hassairi Abdelhamid}$ has defined and
characterized the so-called class of multivariate simple cubic NEF's
which is the natural extension of the class of real cubic NEF's. It
is worth mentioning here that the simple quadratic NEF's are not the
only families which have quadratic variance functions, the Wishart
families on symmetric matrices have also quadratic variance
functions. The classifications of NEF's by the form of the variance
function provide an important tool in the study of distributions. In
fact, in many important cases, the variance function is very simple
and is easier to use than the distribution itself or the Laplace
transform. \n Moreover, the fact that the variance function is
quadratic or cubic, is not only a question of form, but the form
corresponds to some very interesting analytical characteristic
properties. In this respect, let us mention that for the Morris
class of real quadratic NEF's, we have the Meixner characterization
based on some families of orthogonal polynomials which generate
exactly the Morris class (see$\cite{Meixner(1934)}$). Another
characterization due to Feinsilver$\cite{Feinsilver (1986)}$ states
that a certain class of polynomials naturally associated to a NEF is
orthogonal if and only if the family is in the Morris class. This
characterization has been extended to the Casalis class of simple
quadratic NEF's by Labeye-Voisin, and Pommeret $\cite{E.}$.
Concerning the cubic NEF's, Hassairi and Zarai $\cite{Hassairi and
Zarai(2004)}$ introduced a notion of 2-orthogonality for a sequence
of polynomials to give an extended version of the Meixner and
Feinsilver characterization which subsume the Letac-Mora class of
real cubic NEF's. Hassairi and Zarai $\cite{Hassairi and
Zarai(2006)}$ have also introduced a notion of trans-orthogonality
for a sequence of multivariate polynomials to extend their
characterization result to the class of multivariate simple cubic
NEF's. Besides these characterizations based essentially on
different notions of orthogonality of polynomials, it is stated in
Casalis$\cite{Casalis(1996)}$ that the simple quadratic NEF
satisfies a property based on two conjugates families of prior
distributions related to the NEF. For a NEF $F=F(\mu)$, consider the
family of prior distributions $\Pi$ introduced by Diaconis and
Ylvisaker $\cite{Diaconis and Ylvisaker(1979)}$ and defined by
\begin{equation}\label{1}
\Pi=\{\pi_{t,m_{0}}(\textrm{d}\theta)=C_{t,m_{0}}\ \exp t(\langle
m_{0},\theta\rangle-k_{\mu}(\theta))\
{\mathbf{1}}_{\Theta(\mu)}(\theta) d\theta, t>0, m_{0}\in M_{F}\}
\end{equation}
where $C_{t,m_{0}}$ is a normalizing constant. Consider also the
family $\Pi^*$ introduced by Consonni $\textit{et al}$
$\cite{Consonni(2004)}$, see also $\cite{Gutierrez(2003)}$ and
defined by
\begin{equation}\label{2}
\Pi^*=\{\pi^*_{t,m_{0}}, t>0, m_{0}\in M_{F}\}
\end{equation}
where
$$\pi^*_{t,m_{0}}(dm)=C^*_{t,m_{0}}\ \exp t(\
\langle m_{0},\psi_{\mu}(m)\rangle-k_{\mu}(\psi_{\mu}(m))\ )\
{\mathbf{1}}_{M_{F}}(m)dm,$$ and the constant $C^*_{t,m_{0}}$ is a
normalizing constant. Then, when $F(\mu)$ is a Wishart or a simple
quadratic NEF, we have that $k'_{\mu}(\Pi)=\Pi^*$. It is also shown
that this property is equivalent to two other properties expressed
in terms of some differential equations satisfied by the cumulant
function $k_{\mu}$. In the real case, the property characterizes the
Morris class of real quadratic NEF's, that is $k'_{\mu}(\Pi)=\Pi^*$
if and only if the NEF $F$ is in the Morris class. In the present
paper, we extend these results to the class of multivariate simple
cubic NEF's. We construct two families of prior distributions
related to a multivariate NEF, and we show that these families
coincide when the NEF is simple cubic. We then show that this
property is equivalent to the fact that the cumulant function is a
solution of some Monge-Amp\`{e}re equation and also equivalent to
the fact that the variance function satisfies a differential
equation. As a corollary, we obtain three new characterizations of
the Letac-Mora class of real cubic NEF's.
\section{Some equivalent properties}
Throughout this section, we suppose that $F=F(\mu)$ is a NEF on a
linear vector space $E$ with dimension $n$. Besides the family $\Pi$
of prior distributions defined in (1.3), we introduce another family
$\widetilde{\Pi}$ of prior distributions. Let $\beta$ be in $E^{*}$
such that the set
$$\widetilde{\Theta}=\{\theta \in \Theta(\mu); 1+\langle\beta,k'_{\mu}(\theta)\rangle >0\}$$ is nonempty, and denote
$\widetilde{M}=k'_{\mu}(\widetilde{\Theta})$. Consider the family of
prior distributions
\begin{eqnarray}\label{}
\widetilde{\Pi}=\{\ \widetilde{\pi}_{t,m_{0}}\ ;\ t\in\reel_{+}^{*}\
\textrm{,}\ m_{0}\in M_{F}\},
\end{eqnarray}
where $$\widetilde{\pi}_{t,m_{0}}(dm)=\widetilde{C}_{t,m_{0}}\
(1+\langle\beta,m\rangle)^{-n-2}\ \exp t\{\langle
m_{0},\psi_{\mu}(m)\rangle-k_{\mu}(\psi_{\mu}(m))\}\
{\mathbf{1}}_{\widetilde{M}}(m)dm.$$ With these notations, we next
state and prove our first main result.

 \begin{theorem}\label{A}
The three following properties are equivalent

\n (1) There exists $(a,b,c)\in E\times\reel^{2}$ such that for all
$m$ in $M_{F}$,
$$\det V_{F}(m)=\ (1+\langle\beta,m\rangle)^{n+2}\ \exp\{\langle
a,\psi_{\mu}(m)\rangle+bk_{\mu}(\psi_{\mu}(m))+c\}.$$ (2) There
exists $(a,b)\in E\times\reel$ such that for all $m$ in $M_{F}$ and
any  basis, $(e_{i})_{i=1}^{n}$ of $E$, with dual basis
$(e_{i}^{*})_{i=1}^{n},$  we have
\begin{equation}\label{A1}
\displaystyle\sum_{i=1}^{n}[V'_{F}(m)(e_{i})]e_{i}^{*}=\displaystyle\frac{n+2}{1+\langle\beta,m\rangle}V_{F}(m)(\beta)+a+bm.
\end{equation}

\n (3) There exists an open subset $\Omega$ of $\reel_{+}^{*}\times
M_{F}$ such that
$$k'_{\mu}(\Pi)=\widetilde{\Pi}_{\Omega}=\{\widetilde{\pi}_{t,m_{0}};
(t,m_{0})\in\Omega\}.$$
\end{theorem}
Note that (1) may be stated in terms of the cumulant function as
there exists $(a,b,c)\in E\times\reel^{2}$ such that for all
$\theta$ in $\Theta(\mu)$,
$$\det k''_{\mu}(\theta)=\ (1+\langle\beta,k'_{\mu}(\theta)\rangle)^{n+2}\ \exp\{\langle
a,\theta \rangle+bk_{\mu}(\theta)+c\},$$ that is the cumulant
function is solution of some Monge-Amp\`{e}re equation (see $\cite
{Monge-Ampere}$).

\begin{Pff} We will show that $(1)\Leftrightarrow$ (2) and
$(1)\Leftrightarrow(3)$.
 \n

$(1)\Rightarrow (2)$ Suppose that $V_{F}(m)$ satisfies (1), then we
have
$$\log \det V_{F}(m)=(n+2)\log(1+\langle\beta,m\rangle)
+\{\langle\psi_{\mu}(m),a\rangle+b\ k_{\mu}(\psi_{\mu}(m))+c\}.$$
Taking the derivative, we get
$$trace(V_{F}^{-1}(m)V_{F}^{'}(m)(.))=\frac{(n+2)\langle\beta,.\rangle}{1+\langle\beta,m\rangle}
+\langle\psi_{\mu}^{'}(m)(.),a\rangle+b\langle
m,\psi_{\mu}^{'}(m)(.)\rangle,$$ which is equivalent to
$$\displaystyle\sum_{i=1}^{n}[V'_{F}(m)(.)V_{F}^{-1}(m)(e_{i})](e_{i}^{*})
=\frac{(n+2)\langle\beta,.\rangle}{1+\langle\beta,m\rangle}
+\langle\psi_{\mu}^{'}(m)(.),a\rangle+b\langle
m,\psi_{\mu}^{'}(m)(.)\rangle.$$ Replacing $(.)$ by $V_{F}(m)(.)$,
and using the condition of symmetry
$$V'(m)(V(m)(\alpha))(\beta)=V'(m)(V(m)(\beta))(\alpha)\ \ \forall\
\alpha,\beta\ \in\ E^{*}\ \ (see\ \cite{G.}, page\ 103),$$ we obtain
\begin{equation}\label{e4}\displaystyle\sum_{i=1}^{n}[V'_{F}(m)(e_{i})(.)](e_{i}^{*})=\frac{(n+2)\langle\beta,V_{F}(m)(.)\rangle}{1+\langle\beta,m\rangle}
+\langle a,(.)\rangle+b\langle m,(.)\rangle.\end{equation} As
$V_{F}(m)$ is symmetric, we get
$$\displaystyle\sum_{i=1}^{n}[V'_{F}(m)(e_{i})](e_{i}^{*})=\frac{(n+2)}{1+\langle\beta,m\rangle}V_{F}(m)(\beta)
+a+bm.$$

\n \ $(2)\Rightarrow(1)$\  Suppose that (2) holds. Then, we easily
get (\ref{e4}).

\n Replacing, in (\ref{e4}),  (.) by $V_{F}^{-1}(m)(.)$, one obtains
$$\displaystyle\sum_{i=1}^{n}V'_{F}(m)(e_{i})V_{F}^{-1}(m)(.)(e_{i}^{*})=\displaystyle\frac{(n+2)\
\langle\beta,(.)\rangle}{1+\langle\beta,m\rangle}+\langle
a,\psi'_{\mu}(m)(.)\rangle+b\langle m,\psi'_{\mu}(m)(.)\rangle.$$

\n This is equivalent to
$$trace(V_{F}^{-1}(m) V'(m)(.))=\displaystyle\frac{(n+2)\
\langle\beta,(.)\rangle}{1+\langle\beta,m\rangle}+\langle
a,\psi'_{\mu}(m)(.)\rangle+b\langle m,\psi'_{\mu}(m)(.)\rangle.$$

\n Integrating, we deduce that there exists $c$ in $\reel$ such that
$$\log\det(V_{F}(m))=(n+2)\log(1+\langle\beta,m\rangle)
+\{\langle\psi_{\mu}(m),a\rangle+b\ k_{\mu}(\psi_{\mu}(m))+c\},$$
and the result follows.

\n \ $(1)\Rightarrow (3)$ Suppose that (1) holds, and define
$$\Omega=\{(t,m_{0})\in\reel_{+}^{*}\times M_{F}\ ; t>b\ and\ m_{0}\in
(1-\displaystyle\frac{b}{t})M_{F}-\displaystyle\frac{a}{t}\}.$$

\n Take $(t,m_{0})$ in $\Omega$ and denote $\nu$ the image of
$\widetilde{\pi}_{t,m_{0}}$ by $\psi_{\mu}$. Then it is easy to
verify that
$$\nu(d\theta)=\widetilde{C}_{t,m_{0}}\ e^{c}\ \exp\{\langle
tm_{0}+a,\theta\rangle-(t-b)k_{\mu}(\theta)\}\
\mathbf{1}_{\widetilde{\Theta}}(\theta)\textit{d}\theta.$$
\\ Since $(t,m_{0})$ is in $\Omega$, we have that
$t-b>0\ and\ \displaystyle\frac{tm_{0}+a}{t-b}\in M_{F}.$

\n Thus taking $t_{1}=t-b$ and
$m_{1}=\displaystyle\frac{tm_{0}+a}{t-b}$, we obtain that
$$\nu(d\theta)=C_{t_{1},m_{1}}\ \exp\{ t(\langle
m_{1},\theta\rangle-k_{\mu}(\theta))\}\
\mathbf{1}_{\widetilde{\Theta}}(\theta) \textit{d}\theta.$$

\n Hence $\psi_{\mu}(\widetilde{\Pi}_{\Omega})\subset\Pi$, and it
follows that $\widetilde{\Pi}_{\Omega}\subset k'_{\mu}(\Pi).$

\n Conversely, if $\pi_{t,m_{0}}$ is an element of $\Pi$, then its
image $\sigma$ by $k'_{\mu}$ is given by
$$\sigma(\textit{dm})= C_{t,m_{0}}\ e^{-c}\
(1+\langle\beta,m\rangle)^{-n-2}\ \exp\{\langle
tm_{0}-a,\psi_{\mu}(m)\rangle-(t+b)k_{\mu}(\psi_{\mu}(m))\}\
\mathbf{1}_{\widetilde{M}}(\textit{m})\textit{dm}.$$

\n Taking $t_{1}=t+b$ and $m_{1}=\displaystyle\frac{tm_{0}-a}{t+b}$.
Then $(t_{1},m_{1})$ is in $\Omega$, and we have
$$\sigma(\textit{dm})=\widetilde{C}_{t_{1},m_{1}}\
(1+\langle\beta,m\rangle)^{-n-2}\ \exp t_{1}\{\langle
m_{1},\psi_{\mu}(m)\rangle-k_{\mu}(\psi_{\mu}(m))\}\mathbf{1}_{\widetilde{M}}(\textit{m})\textit{dm}\
,$$ which is an element of $\widetilde{\Pi}_{\Omega}.$ \n \

$(3)\Rightarrow (1)$ Suppose that
$k'_{\mu}(\Pi)=\widetilde{\Pi}_{\Omega}$. Then, for an element
$\pi_{t,m_{0}}$ of $\Pi$, we have on the one hand,
$$k'_{\mu}(\pi_{t,m_{0}})(\textit{dm})=(\det V_{F}(m))^{-1}\ C_{t,m_{0}}\
\exp t\{\langle m_{0},\psi_{\mu}(m)\rangle-k_{\mu}(\psi_{\mu}(m))\}\
\mathbf{1}_{\widetilde{M}}(\textit{m})\textit{dm}.$$ On the other
hand, since $k'_{\mu}(\pi_{t,m_{0}})$ is in
$\widetilde{\Pi}_{\Omega}$, there exists $(t_{1},m_{1})$ in $\Omega$
such that
$$k'_{\mu}(\pi_{t,m_{0}})(\textit{dm})=\widetilde{C}_{t_{1},m_{1}}\
(1+\langle\beta,m\rangle)^{-n-2}\ \exp t_{1}\{\langle
m_{1},\psi_{\mu}(m)\rangle-k_{\mu}(\psi_{\mu}(m))\}\
\mathbf{1}_{\widetilde{M}}(\textit{m})\textit{dm}.$$ Comparing these
two expressions of $k'_{\mu}(\pi_{t,m_{0}})$ gives $$\ \det
V(m)=(1+\langle\beta,m\rangle)^{n+2}\ \exp\{\langle
a,\psi_{\mu}(m)\rangle+bk_{\mu}(\psi_{\mu}(m))+c\},$$ where
$a=tm_{0}-t_{1}m_{1}$, $b=t_{1}-t$, and $c=\log
(\displaystyle\frac{c_{t,m_{0}}}{\widetilde{c}_{t_{1},m_{1}}})$.

\end{Pff}

\section{Characterizations of the Letac-Mora class of real cubic NEFs}
In this section, we prove that a multivariate simple cubic NEF
satisfies the properties in Theorem (\ref{A}), and that the real
versions of these properties characterize the real cubic NEFs.
Recall that a simple cubic NEF is obtained form a simple quadratic
NEF by the so-called action of the linear group $GL(\reel\times E)$
on the NEFs of $E$. For more details, we refer the reader to
$\cite{Hassairi Abdelhamid}$, where a complete description of this
class is given. This action is in fact an extension of the way in
which the Letac-Mora class of real cubic NEFs is obtained from the
Morris class of real quadratic NEF's. For our purposes here, we need
only to mention that, up to affine transformations and power of
convolution, a simple cubic variance function is of the form
\begin{equation}\label{e1}V(m)=(1+\langle\beta,m\rangle)\
(I+m\otimes\beta)\
V_{1}(\displaystyle\frac{m}{1+\langle\beta,m\rangle})\
(I+\beta\otimes m),\end{equation} where $V_{1}$ is the variance
function of a simple quadratic NEF $F_{1}$, and $m$ is in
$(M_{F_{1}})_{\beta}$, where
$$(M_{F_{1}})_{\beta}=\{m\in M_{F_{1}}; 1+\langle\beta,m\rangle>0\ and\ \displaystyle\frac{m}{1+\langle\beta,m\rangle}\in M_{F_{1}}\}.$$ The relation (\ref{e1}) is invertible and conversely, we have
\begin{equation}\label{e2}V_{1}(M)=(1-\langle\beta,M\rangle)\ (I-M\otimes\beta)\
V(\displaystyle\frac{M}{1-\langle\beta,M\rangle})\ (I-\beta\otimes
M),\end{equation} where $M$ is in $(M_{F})_{-\beta}$.

\n We also mention that the relation between a simple cubic NEF
$F(\mu)$ and a simple quadratic NEF $F(\nu)$ may also be expressed
in terms of the cumulant functions by

\n \

\begin{equation} \label{e6}\left\{\begin{array}{ccc}k_{\mu}(\lambda)=k_{\nu}(\theta)-k_{0}\\\
\\ \lambda=-\beta
k_{\nu}(\theta)+\theta-\lambda_{0}\end{array}\ \
\right.\end{equation}

\n or equivalently by

\begin{equation} \left\{\begin{array}{ccc}k_{\nu}(\theta)=k_{\mu}(\lambda)-k_{1}\\\
\\ \theta=\beta
k_{\mu}(\lambda)+\lambda-\theta_{1}\end{array}\ \
\right.\end{equation}

\n where $(k_{0},\lambda_{0})$ and $(k_{1},\theta_{1})$ are
constants in $\reel\times E.$
 \n Note that if $\beta=0$ in
(\ref{e1}) we obtain the simple quadratic class. Then for more
accuracy we exclude this case and we keep only $\beta$ in
$E^{*}\setminus\{0\}$.

\n We now prove that the multivariate simple cubic NEF's satisfy the
properties in Theorem \ref{A}.

\begin{prop}$\label{A3}$
Let $F=F(\mu)$ be a simple cubic NEF on $E$, then there exists
$(a,b,c)$ in $E^{*}\times\reel^{2}$ such that
$$\det(V_{F}(m))=(1+\langle\beta,m\rangle)^{n+2}\ \exp\{\langle\psi_{\mu}(m),a\rangle+b\
k_{\mu}(\psi_{\mu}(m))+c\}.$$\end{prop}

\begin{Pff} Given that the family $F$ is simple cubic, then there exist $\beta$ in $E^{*}$ and
$F_{1}=F(\nu)$ a simple quadratic NEF such that
$$V_{F}(m)=(1+\langle\beta,m\rangle)\
(I+m\otimes\beta)\
V_{F_{1}}(\displaystyle\frac{m}{1+\langle\beta,m\rangle})\
(I+\beta\otimes m),$$ see  (\ref{e1}). As
$\det(I+m\otimes\beta)=1+\langle\beta,m\rangle$, we obtain

$$\det(V_{F}(m))=(1+\langle\beta,m\rangle)^{n+2}\
\det(V_{F_{1}}(\displaystyle\frac{m}{1+\langle\beta,m\rangle})).$$

\n We now use the fact for a simple quadratic NEF $F_{1}$ (see
$\cite {Casalis(1996)})$, there exist $a'$ in $E^{*}$ and $b',c'$ in
$\reel$ such that, for all $M$ in $M_{F_{1}}$,
$$\det(V_{F_{1}}(M))=\exp\{\langle a',\psi_{\nu}(M)\rangle+b'
k_{\nu}(\psi_{\nu}(M))+c'\}.$$ It follows that
$$\det(V_{F}(m))=(1+\langle\beta,m\rangle)^{n+2}\
\exp\{\langle\psi_{\nu}(\displaystyle\frac{m}{1+\langle\beta,m\rangle}),a'\rangle+b'\
k_{\nu}(\psi_{\nu}(\displaystyle\frac{m}{1+\langle\beta,m\rangle}))+c'\}.$$

\n From (\ref{e6}), putting $\lambda=\psi_{\mu}(m)$ and
$\theta=\psi_{\nu}(\ \displaystyle\frac{m}{1+\langle\beta,m\rangle}\
)$, we get
$$k_{\nu}(\psi_{\nu}(\ \displaystyle\frac{m}{1+\langle\beta,m\rangle}\ ))=k_{\mu}(\psi_{\mu}(m))+k_{0},$$

$$\psi_{\nu}(\ \displaystyle\frac{m}{1+\langle\beta,m\rangle}\ )=\psi_{\mu}(m)+\beta
k_{\mu}(\psi_{\mu}(m))+\beta k_{0}+\lambda_{0}.$$ Then
\begin{eqnarray*}
\det V_{F}(m)&=&(1+\langle\beta,m\rangle)^{n+2}
\exp\{\langle\psi_{\mu}(m),a'\rangle+(b'+\langle
a',\beta\rangle)k_{\mu}(\psi_{\mu}(m))\\& &+(b'k_{0}+\langle
a',\beta k_{0}+\lambda_{0}\rangle+c') \}.\end{eqnarray*}
 \n Setting $a=a'$, $b=b'+\langle a',\beta\rangle$ and $c=b'k_{0}+\langle
a',\beta k_{0}+\lambda_{0}\rangle+c'$, we obtain the desired result.
\end{Pff}

As the Letac-Mora class of real cubic NEFs is nothing but the simple
cubic class, when the dimension $n$ is equal to 1, this class
satisfies the real version of the properties in Theorem \ref{A}. We
will show that, in this case, these properties are characteristic.

\begin{theorem}
 Let $F=F(\mu)$ be a NEF on the real line, then $F$ is cubic if
and only if $k'_{\mu}(\Pi)=\widetilde{\Pi}$.

\end{theorem}

\begin{Pff}
Suppose that $k'_{\mu}(\Pi)=\widetilde{\Pi}$. Then according to
Theorem(\ref{A}), the variance function $V_{F}(m)$ satisfies the
differential equation
$$(1+\beta m) V_{F}'(m)-3\beta\ V_{F}(m)=(a+bm)(1+\beta m).$$
Solving this differential equation by standard methods gives

$$V_{F}(m)=\lambda\ (1+\beta m)^{3}-\displaystyle\frac{b}{\beta^{2}}\ (1+\beta m)^{2}+\displaystyle\frac{b-\beta a}{2\beta^{2}}\ (1+\beta m),$$ which is a polynomial of degree less then or equal to
3.
\end{Pff}

\end{document}